\documentclass[a4paper, 10pt, conference]{ieeeconf}      

\IEEEoverridecommandlockouts                              

\usepackage{graphics} 
\usepackage{epsfig} 
\usepackage{mathptmx} 
\usepackage{times} 
\usepackage{amsmath} 
\usepackage{amssymb}  
\usepackage{comment}
\usepackage{mathrsfs}
\usepackage{eucal}
\usepackage{cuted}
\usepackage{float}
\let\labelindent\relax
\usepackage{nicefrac}
\usepackage{enumitem}
\usepackage{verbatim}
\usepackage{siunitx}
\usepackage{textcomp}
\usepackage{tikz}
\usepackage{savesym}
\usepackage{amsfonts}
\usepackage{mathtools}
    \DeclarePairedDelimiter\norm{\lVert}{\rVert} 
\usepackage{graphicx}
\usepackage{xcolor}
\usepackage{booktabs}
\usepackage{tabularx}
\usepackage{array}

\newtheorem{conjecture}{Open question}

\usetikzlibrary{shapes,arrows,positioning}
\usepackage[normalem]{ulem}

\usepackage{url}

\newcommand{\Ln}{\mathrm{Ln}}
\newcommand{\Exp}{\mathrm{Exp}}
\newcolumntype{Y}{>{\raggedright\arraybackslash}X}
\newcolumntype{P}[1]{>{\raggedright\arraybackslash}p{#1}}

\title{\LARGE \bf
Port-Hamiltonian Systems Beyond Joules:\\ A Review of Unexpected Application Domains}

\author{Mauricio Mu\~noz-Arias$^{1}$
\thanks{$^{1}$Mauricio Mu\~noz-Arias is with the Engineering and Technology Institute Groningen (ENTEG), Faculty of Science and Engineering, University of Groningen, Groningen, The Netherlands.
        {\tt\small m.munoz.arias@rug.nl}}%
}

\begin{document}

\maketitle
\thispagestyle{empty}
\pagestyle{empty}

\begin{abstract}
The port-Hamiltonian framework was built for mechanical and electrical systems, where the Hamiltonian is stored energy and the ports carry power. Over the past fifteen years it has spread to chemistry, cell biology, fusion plasmas, musical instruments, district heating, electricity markets, macroeconomics, predator-prey ecology, optimization algorithms and neural-network training. This review maps that spread. We classify each application by what its Hamiltonian really is: physical energy, a thermodynamic potential such as free energy or exergy, or an abstract storage function with no energetic meaning. The review keeps only sources whose bibliographic record we checked against Crossref. The picture that emerges is a continuum. Where a genuine thermodynamic potential exists, the theory is deep and constructive. Where the energy is borrowed, the framework still pays off through compositional interconnection and dissipation inequalities, but the work is younger and scattered. We close with a set of deliberately speculative research directions, including open port-Hamiltonian models of epidemics, opinion dynamics and supply chains, and learned port-Hamiltonian digital twins for structural health monitoring.
\end{abstract}

\section{INTRODUCTION}

Passivity-based control rests on one physical intuition.
Systems that store and dissipate energy settle where their energy is lowest, so a controller that reshapes the energy and adds damping can place the equilibrium where we want it \cite{ortega2001energy, ortega2002idapbc}.
The port-Hamiltonian (pH) framework turns that intuition into a modelling language.
A system is described by an energy function, a power-preserving interconnection structure, a dissipation structure and a set of power ports \cite{vds2017l2gain, duindam2009, vdsJeltsema2014}.

The language was born in mechanics and circuits.
Its reach is now much wider.
Chemical reaction networks have a pH form in which the Hamiltonian is Gibbs free energy \cite{vdsRaoJay2013}.
Tokamak plasmas have been controlled with interconnection and damping assignment on a pH model \cite{vu2016cep}.
Analog audio circuits and electric pianos are simulated with schemes that keep the power balance exact in discrete time \cite{falaize2016, falaize2017rhodes}.
A primal-dual pricing algorithm has been written as a pH system and plugged into a power-grid model \cite{stegink2017tac}.
Even macroeconomic flows \cite{macchelli2014econ}, predator-prey populations \cite{massaroli2020lv} and the training of neural networks \cite{massaroli2019cdc} have been given pH structure.
Our own work followed part of this path, from force and visual-servo control of robot manipulators \cite{munozarias2013cdc, munozarias2015ibvs} to wave energy converters \cite{barradas2017og, almuzakki2018cummins}, pumped-hydro storage \cite{phillips2019pumped, esquivel2022pumpturbine} and drag-free satellites for gravitational-wave detection \cite{ma2024gw}.

These applications differ in one essential respect.
In some of them the Hamiltonian is measured in joules.
In others it is a thermodynamic potential that has physical status but is not total energy.
In the most surprising ones it is a Lyapunov-like function with no energetic meaning at all.
The question this review asks is simple: how far has the framework been stretched, and what does it still deliver once the joules are gone?

\textit{Contributions.}
First, we propose a three-band classification of pH applications by the nature of the Hamiltonian, and use it to organise the literature.
Second, we compile a verified repository of sources across fourteen groups of non-classical domains, with every DOI checked against its Crossref record.
Third, we propose concrete open questions, some of them deliberately speculative, where pH theory looks natural but has not yet been applied.

The paper is organised as follows.
Section~\ref{sec:Preliminaries} recalls the pH framework and defines the three bands.
Section~\ref{sec:Method} describes how the review was built.
Sections~\ref{sec:Thermo}, \ref{sec:Physical} and \ref{sec:Abstract} review the thermodynamic, physical and abstract bands.
Section~\ref{sec:Repository} gives the master table and a teaching shortlist.
Section~\ref{sec:Crazy} collects open questions, and Section~\ref{sec:Concluding_Remarks} concludes.

\section{PRELIMINARIES}
\label{sec:Preliminaries}

\subsection{The port-Hamiltonian framework}
An input-state-output pH system is given by \cite{vds2017l2gain, vdsJeltsema2014}
\begin{equation}
\Sigma_{PH}
\left\{
\begin{aligned}
    \dot{x}
    &=
    \left[ J(x)-R(x) \right]
    \nabla_{x}H(x)
    +
    g(x)\,u,
    \\
    y &=g(x)^\top
    \nabla_{x}H(x),
\end{aligned}
\right.
\label{eq:sigma_total}
\end{equation}
with state $x \in \mathbb{R}^{n}$, input $u \in \mathbb{R}^{m}$ and output $y \in \mathbb{R}^{m}$.
The interconnection matrix is skew-symmetric, $J(x) = -J(x)^\top$, and the dissipation matrix is positive semi-definite, $R(x) = R(x)^\top \succeq 0$.
Differentiating the Hamiltonian $H(x)$ along \eqref{eq:sigma_total} gives the power balance
\begin{equation}
\dot{H} =
- \nabla_{x}^\top H \, R \, \nabla_{x}H
+
y^\top u
\;\leq\;
y^\top u .
\label{eq:power_balance}
\end{equation}
The system is therefore passive with respect to the port pair $(u,y)$ whenever $H$ is bounded from below.
Passivity-based control exploits \eqref{eq:power_balance} in two steps.
Energy shaping replaces $H$ by a closed-loop function $H_d$ with a strict minimum at the desired equilibrium.
Damping injection then adds dissipation so that $H_d$ decreases along trajectories \cite{ortega2001energy, ortega2002idapbc}.

Three structural facts matter for what follows.
The power-conserving interconnection of pH systems is again a pH system, with the sum of the Hamiltonians as storage.
This compositional property does not depend on what $H$ means.
Distributed-parameter systems fit the same mould through Stokes-Dirac structures and boundary ports \cite{rashad2020review}.
Systems with algebraic constraints fit through pH descriptor formulations \cite{mehrmannUnger2023}.

\subsection{Three meanings of the Hamiltonian}
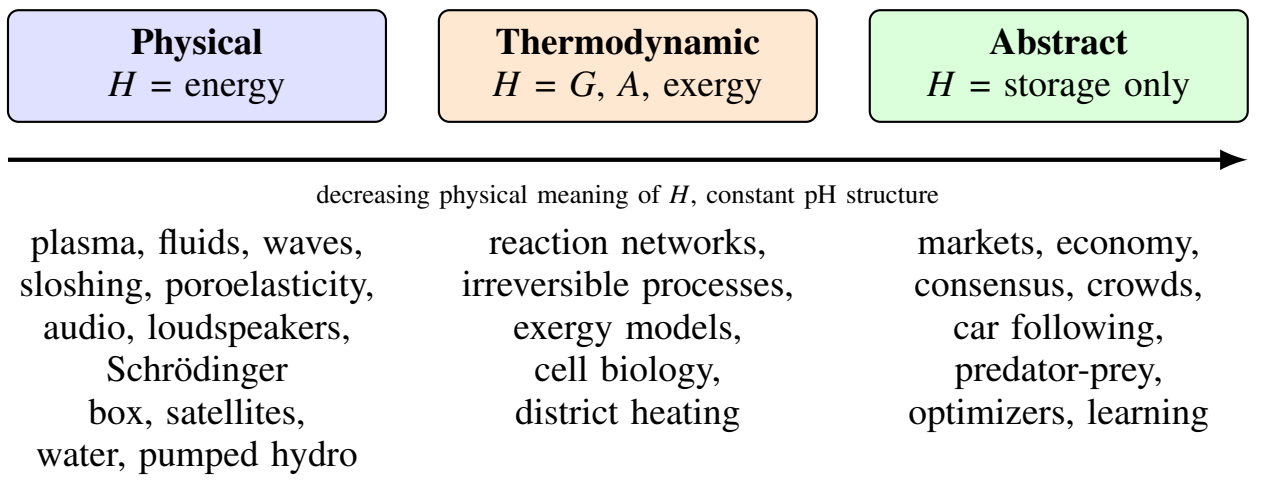
\begin{figure}[t]
\centering
\begin{tikzpicture}[
  font=\scriptsize,
  band/.style={draw, rounded corners=2pt, minimum width=2.5cm, minimum height=0.55cm, align=center},
  dom/.style={align=center, text width=2.5cm}
]
\node[band, fill=blue!12]   (P) at (0,0)    {\textbf{Physical}\\ $H$ = energy};
\node[band, fill=orange!18] (T) at (2.85,0) {\textbf{Thermodynamic}\\ $H$ = $G$, $A$, exergy};
\node[band, fill=green!14]  (A) at (5.7,0)  {\textbf{Abstract}\\ $H$ = storage only};
\draw[-latex, thick] (-1.25,-0.62) -- (6.95,-0.62);
\node[font=\tiny, anchor=north] at (2.85,-0.66) {decreasing physical meaning of $H$, constant pH structure};
\node[dom, anchor=north] at (0,-0.95)    {plasma, fluids, waves,\\ sloshing, poroelasticity,\\ audio, loudspeakers,\\ Schr\"odinger box, satellites,\\ water, pumped hydro};
\node[dom, anchor=north] at (2.85,-0.95) {reaction networks,\\ irreversible processes,\\ exergy models,\\ cell biology,\\ district heating};
\node[dom, anchor=north] at (5.7,-0.95)  {markets, economy,\\ consensus, crowds,\\ car following,\\ predator-prey,\\ optimizers, learning};
\end{tikzpicture}
\caption{The three meanings of the Hamiltonian in the reviewed literature. The interconnection, dissipation and port structure of \eqref{eq:sigma_total} is the invariant; the interpretation of $H$ changes from left to right.}
\label{fig:continuum}
\end{figure}

We sort applications into three bands according to what $H$ represents.
\begin{itemize}[leftmargin=*]
\item \textbf{Physical (P).} $H$ is stored kinetic, potential, electromagnetic, acoustic, hydraulic or internal energy. The ports carry power in watts.
\item \textbf{Thermodynamic (T).} $H$ is a thermodynamic potential that is not total energy: Gibbs free energy, availability, exergy, or an entropy-based function. It has physical status and imposes constitutive restrictions, but energy conservation and entropy production must be kept apart.
\item \textbf{Abstract (A).} $H$ is a first integral, a convex objective, a disagreement measure or another Lyapunov-like storage function. The ports carry a bilinear supply rate that need not be power.
\end{itemize}
Fig.~\ref{fig:continuum} places the reviewed domains on this continuum.
The bands are not watertight.
A power grid coupled to a market controller is physical on one side of the interconnection and abstract on the other \cite{stegink2017tac}.
A pH neural network trained on mechanical data learns a function that is meant to approximate energy, but nothing forces it to \cite{desai2021}.
Where a source sits on a boundary, the master table gives both labels.

\section{HOW THIS REVIEW WAS BUILT}
\label{sec:Method}

\subsection{Search and screening}
Candidate sources were gathered with AI-assisted literature searches and screened by hand.
We looked for peer-reviewed applications of pH theory outside rigid-body mechanics, robotics, circuits, hydraulics and pneumatics, and recorded for each source its Hamiltonian, its port variables and its energy band.

\subsection{Verification and inclusion rules}
Every source in the bibliography was resolved against its Crossref record in October and November 2026, checking authors, title, venue, volume, pages and year.
The single exception is a PMLR paper, which has no DOI and was checked on the publisher page \cite{neary2023}.
Three author lists were corrected in this step.
The non-isothermal reaction network paper is by Wang, Maschke and van der Schaft, not by van der Schaft alone \cite{wang2018jmc}.
The first author of the cardiac bond-graph paper is Pan \cite{pan2018cardiac}.
The ideal-fluid paper has four authors, not five \cite{rashad2021fluid}.
We excluded sources we could not confirm and recent preprints that have not yet been peer reviewed, even when their topic was attractive.
These include 2025 and 2026 preprints on Hopfield networks, quantum pH neural networks, cortical and cellular dynamics, and Hamiltonian reformulations of epidemic models.
Bond-graph papers are included where they are the clear counterpart of a pH formulation, as in systems biology.

\section{THE THERMODYNAMIC BAND}
\label{sec:Thermo}

\subsection{Chemical reaction networks}
This is the most complete non-mechanical branch of pH theory.
Van der Schaft and Maschke first gave open reaction networks a pH form \cite{vdsMaschke2010crn}.
Van der Schaft, Rao and Jayawardhana then wrote detailed-balanced mass-action kinetics as \cite{vdsRaoJay2013, vdsRaoJay2016}
\begin{equation}
\dot{x} = -Z \mathcal{L} \, \Exp\!\left( Z^\top \Ln \frac{x}{x^*} \right),
\label{eq:crn}
\end{equation}
where $x$ holds the species concentrations, $x^*$ is a thermodynamic equilibrium, $Z$ is the complex composition matrix and $\mathcal{L}$ is a balanced weighted Laplacian of the complex graph.
The vector $\Ln(x/x^*)$ is the vector of chemical potentials, and it is the gradient of the dimensionless Gibbs free energy
\begin{equation}
G(x) = \sum_{i} \left( x_i \ln \frac{x_i}{x_i^*} - x_i + x_i^* \right).
\label{eq:gibbs}
\end{equation}
Along \eqref{eq:crn}, $\dot{G} \leq 0$, with equality only on the equilibrium set.
Classical results of reaction network theory, such as the uniqueness of equilibria in each stoichiometric class, then follow from the dissipation inequality and the Laplacian structure.
The function \eqref{eq:gibbs} will reappear in Section~\ref{sec:Crazy}.

Šešlija, van der Schaft and Scherpen extended the formulation to reaction-diffusion systems \cite{seslija2014}.
They started from the distributed pH model and derived a compartmental discretisation that keeps the thermodynamic structure.
Wang, Maschke and van der Schaft removed the isothermal assumption by adding internal energy as a state \cite{wang2016fosbe, wang2018jmc}.
Their example is a genetic protein-synthesis circuit with cell-to-cell communication, which shows how directly the theory reaches biology.

\subsection{Irreversible pH systems}
Standard pH models put all irreversibility into $R(x)$.
In thermodynamics this is not enough.
Heat conduction between two bodies conserves total energy exactly yet produces entropy, so a term that subtracts energy cannot describe it.
Ramirez, Maschke and Sbarbaro introduced irreversible pH systems to resolve this \cite{ramirez2013ces}.
In their simplest form,
\begin{equation}
\dot{x} = \gamma(x)\,\{S,U\}_J \, J \, \nabla U(x) + g(x,\nabla U)\,u,
\label{eq:iphs}
\end{equation}
where $U$ is internal energy, $S$ is entropy, $\{S,U\}_J = \nabla S^\top J \nabla U$ is a pseudo-bracket and $\gamma(x) \geq 0$ encodes the irreversible constitutive law.
Energy is conserved by skew-symmetry and entropy production is non-negative by construction.
The continuous stirred tank reactor was the driving example.

The line has grown steadily.
It covers multi-energy systems \cite{ramirez2013ejc}, electromechanical actuators with hysteresis \cite{ramirez2018hyst}, boundary-controlled distributed processes \cite{ramirez2022ces} and modulated interconnection of irreversible subsystems \cite{ramirez2024aut}.
An overview appeared in 2022 \cite{ramirez2022entropy}.

The control result is the part most worth teaching.
For an isolated thermodynamic system total energy is conserved, so it cannot serve as a Lyapunov function.
Ramirez et al. used the thermodynamic availability function, a potential measured relative to a reference state, as the closed-loop storage function \cite{ramirez2013ifac, ramirez2016aut}.
Lohmayer, Kotyczka and Leyendecker pushed the same idea further and argued that exergy,
\begin{equation}
E = U - T_0 S + p_0 V - \textstyle\sum_k \mu_{0,k} N_k ,
\label{eq:exergy}
\end{equation}
is often the natural Hamiltonian of a thermodynamic pH model, because it measures the useful work extractable relative to an environment at $(T_0, p_0, \mu_0)$ \cite{lohmayer2021}.
Delvenne and Sandberg connected pH passivity with finite-time thermodynamic bounds \cite{delvenne2014}.
Van der Schaft's tutorial recasts classical thermodynamics in the language of systems and control and is an excellent entry point for students \cite{vds2021csm}.

\subsection{Cell biology and physiology}
The richest biological line uses bond graphs, the graphical counterpart of pH systems.
Gawthrop, Crampin and co-workers model biochemical networks with chemical potential as effort and molar flow as flow, so that stoichiometry and thermodynamics cannot silently contradict each other \cite{gawthrop2014, gawthrop2015, gawthrop2016iet, gawthrop2021mod}.
The models cover biochemical cycles, hierarchical pathways, the cardiac action potential \cite{pan2018cardiac} and advection of species by a carrier fluid \cite{gawthrop2022adv}.
Jeltsema and van der Schaft gave memristive elements a pH home, motivated in part by the ion channels of the Hodgkin-Huxley neuron \cite{jeltsema2010}.
We found no peer-reviewed pH model of whole-body cardiovascular regulation.

\section{THE PHYSICAL BAND IN UNUSUAL PLACES}
\label{sec:Physical}

\subsection{Fusion plasmas}
Tokamak plasma control is perhaps the most striking physical application.
Vu, Lefèvre and Maschke derived a pH model of tokamak thermo-magneto-hydrodynamics that couples Maxwell's equations with mass, momentum and entropy balances and irreversible entropy production \cite{vu2016mcmds}.
On a reduced model of the coupled magnetic flux and heat diffusion, they designed an IDA-PBC controller \cite{vu2014ifac}.
The safety-factor profile controller was then tested on the TCV tokamak \cite{vu2016cep}.
This is one of the few cases in this review where pH control has been run on a large experimental facility.

\subsection{Space instruments}
Deep space is another extreme setting.
Drag-free satellites for gravitational-wave detection must keep their test masses almost perfectly free of disturbance, and their dynamics are strongly nonlinear.
Ma et al. modelled the coupled spacecraft and test-mass dynamics in pH form and designed an energy-based controller with integral action that rejects unknown disturbances and gives input-to-state stability \cite{ma2024gw}.
Linearised designs, the usual choice in this field, offer no such guarantee away from the operating point.

\subsection{Fluids, waves and poroelasticity}
The distributed pH literature is surveyed in \cite{rashad2020review}.
Ideal fluid flow on Riemannian manifolds has a rigorous Stokes-Dirac formulation \cite{rashad2021fluid}.
Reactive Navier-Stokes equations \cite{altmannSchulze2017} and two-phase flow models \cite{bansal2021} have pH forms, although the general drift-flux model may not be energy-consistent \cite{bansal2021}.
Liquid sloshing in moving containers has been modelled as a pH system and coupled to a flexible structure \cite{cardoso2017}.
Biot poroelasticity and multiple-network brain models have pH descriptor formulations that help discretisation and model reduction \cite{altmann2021poro}.
Ocean wave energy is a further case.
Cummins' equation, which describes the radiation forces on a floating body, has been given a pH form, so that arrays of floaters and their power take-off can be interconnected and checked against the WEC-Sim toolbox \cite{almuzakki2018cummins}.
It builds on a modular pH model of the power take-off of the Ocean Grazer wave energy converter \cite{barradas2017og}.

\subsection{Audio and musical acoustics}
Falaize and Hélie convert an analog audio schematic into a pH system automatically: the netlist becomes a graph, the graph becomes a Dirac structure, and the storage and dissipative elements are read from the components \cite{falaize2016}.
The model is simulated with a discrete gradient, so the power balance \eqref{eq:power_balance} holds exactly in discrete time.
Passivity then guarantees numerical stability, which matters for real-time audio.
A second-order explicit variant followed \cite{lopes2015}, and the method was applied to a Rhodes electric piano \cite{falaize2017rhodes} and to the nonlinear electrodynamic loudspeaker \cite{falaize2020speaker}.
Here the energy is entirely physical.
What surprises is the payoff: passivity becomes a numerical design principle for sound synthesis.

\subsection{Quantum and infrastructure systems}
Verrier, Haine and Matignon wrote the Schrödinger equation on a bounded domain as a boundary-controlled pH system and simulated Rabi oscillations in a controlled quantum box \cite{verrier2023}.
Open quantum systems with Lindblad dissipation remain untouched.

Infrastructure networks are a growing application area.
Hauschild et al. built a hierarchy of pH models for district heating, from PDE networks to semi-discrete models \cite{hauschild2020}.
Strehle et al. used shifted passivity to prove stability of hydraulic equilibria in fourth-generation district heating with distributed pumps \cite{strehle2021}.
Machado, Cucuzzella and Scherpen proved passivity properties for a meshed network with several producers, consumers and storage tanks \cite{machado2022}.
Krishna and Schiffer coupled electrical and thermal energy in a microgrid \cite{krishna2021}.
Perryman, Taylor and Karney used an incremental pH form to design decentralised pump controllers for water distribution networks \cite{perryman2022}.
Small pumped-hydro plants join three domains in one storage function.
A reversible pump-turbine driven by an asynchronous machine has a single Hamiltonian over electrical, mechanical and hydraulic energy, and energy shaping controls it in both pumping and generating modes \cite{esquivel2022pumpturbine}.
A solar-powered version of such a plant is modelled in \cite{phillips2019pumped}.

\section{THE ABSTRACT BAND}
\label{sec:Abstract}

\subsection{Networks, consensus, crowds and traffic}
Van der Schaft and Maschke showed that pH systems on graphs cover mass-spring networks, consensus, clustering and distribution networks with one template \cite{vdsMaschke2013graphs}.
For consensus, the Hamiltonian is a quadratic disagreement function and the efforts need not be forces or voltages.
The paper itself notes that many network models with a non-physical background fit the same form.
Wei and van der Schaft used it to design proportional-integral flow controllers for distribution networks with flow constraints and unknown demand \cite{wei2013scl}, a result that reads naturally as inventory control.

Tordeux and Totzeck described pedestrian crowds with multi-scale pH systems \cite{tordeux2023}.
Rüdiger, Tordeux and Ugurcan formulated stochastic car-following models as input-state-output pH systems \cite{rudiger2024}.
Their Hamiltonian is built from headway potentials and relative speeds, and it uses the distance to the following vehicle as an extra degree of freedom.
The added pH term damps stop-and-go waves without lowering the equilibrium speed.
Ackermann et al. added control ports for stabilisation under noise \cite{ackermann2024}.
For autonomous vehicles this is a direct and teachable example.

\subsection{Markets and the economy}
Stegink, De Persis and van der Schaft showed that the gradient method for social-welfare maximisation has a pH form \cite{stegink2015ifac}.
They then interconnected this market controller with a physical pH model of the power grid and proved stability of the closed loop from passivity \cite{stegink2017tac}.
The storage is physical on the grid side and abstract on the market side, and the interconnection does not care.
Van der Schaft presents this closed loop as a true cyber-physical system \cite{vds2020annrev}.

Macchelli's macroeconomic paper goes further from physics than any other source in this review \cite{macchelli2014econ}.
Production, inventory, supply and demand are modelled with effort-flow pairs whose product is a monetary flow rather than power.
Nothing in the model is mechanical, yet the port-based interconnection and the structural analysis carry over.

\subsection{Optimization algorithms}
Optimization is becoming a genuine pH domain in its own right.
Consider the primal-dual gradient flow for $\min_z f(z)$ subject to $Az = b$, with $f$ convex:
\begin{equation}
\dot{z} = -\nabla f(z) - A^\top \lambda, \qquad \dot{\lambda} = A z - b .
\label{eq:pd}
\end{equation}
With $H = \tfrac{1}{2}\norm{z - z^*}^2 + \tfrac{1}{2}\norm{\lambda - \lambda^*}^2$ around a saddle point $(z^*, \lambda^*)$, the coupling through $A$ is skew-symmetric and the gradient term acts as a monotone resistor, so that $\dot{H} = -(z-z^*)^\top(\nabla f(z) - \nabla f(z^*)) \leq 0$.
Camlibel and van der Schaft made this precise.
Incrementally pH systems built from maximal monotone relations are equilibrium-independent passive, and their equilibria can be computed by convex optimisation \cite{camlibel2023}.
Gernandt and Schaller extended the primal-dual pH structure to infinite-dimensional optimal control and interconnected the optimiser with a pH plant \cite{gernandt2025}.
The optimiser stops being a software layer outside the control loop and becomes a dynamical component inside it.

\subsection{Population dynamics}
The Lotka-Volterra predator-prey model has no stored energy.
It does have a classical first integral.
With prey $x_1$, predators $x_2$ and $\dot{x}_1 = a x_1 - b x_1 x_2$, $\dot{x}_2 = -c x_2 + d x_1 x_2$, the change of coordinates $z = \Ln(x)$ gives
\begin{equation}
\dot{z} =
\begin{bmatrix} 0 & -1 \\ 1 & 0 \end{bmatrix}
\nabla H(z),
\quad
H(z) = d e^{z_1} - c z_1 + b e^{z_2} - a z_2 ,
\label{eq:lv}
\end{equation}
a lossless Hamiltonian system.
Massaroli et al. exploited this kind of structure to design passivity-based harvesting and stocking laws that stabilise coexistence without driving either species extinct \cite{massaroli2020lv}.
Their storage has nothing to do with ecological energy or biomass.
It is perhaps the clearest example in the literature of pH reasoning surviving after the physical meaning has been removed entirely.

\subsection{Machine learning}
Two philosophies coexist.
The first uses pH structure as a prior for learning models from data.
Desai et al. learned the Hamiltonian, dissipation and time-dependent forcing with pH neural networks \cite{desai2021}.
Beckers et al. learned a probability distribution over Hamiltonians with Gaussian processes while keeping passivity \cite{beckers2022cdc}.
Neary and Topcu learned subsystems separately and interconnected them afterwards, keeping the pH structure \cite{neary2023}.
Moradi et al. added output-error noise models so that pH networks can be trained on noisy measured data \cite{moradi2026}.

The second philosophy is stranger.
Massaroli et al. treated the weights of a neural network as the state of a pH system and used passivity to analyse the convergence of training \cite{massaroli2019cdc}.
Here the Hamiltonian is an algorithmic convergence measure.
Learning itself becomes a passive dynamical system.

\section{VERIFIED REPOSITORY}
\label{sec:Repository}

Table~\ref{tab:master} lists the core application sources, sorted by domain and then by year.
The band column uses the labels of Section~\ref{sec:Preliminaries}.
The use column distinguishes modelling and analysis (M), control (C), structure-preserving simulation (S) and learning (L).
Surveys and textbooks \cite{vds2017l2gain, duindam2009, vdsJeltsema2014, vds2020annrev, rashad2020review, mehrmannUnger2023} are not repeated in the table.

\begin{table*}[t]
\centering
\caption{Verified repository of non-classical port-Hamiltonian applications. Band: P physical, T thermodynamic, A abstract. Use: M modelling or analysis, C control, S structure-preserving simulation, L learning.}
\label{tab:master}
\scriptsize
\renewcommand{\arraystretch}{0.92}
\begin{tabularx}{\textwidth}{@{}P{2.6cm}P{0.55cm}P{0.9cm}YP{0.5cm}P{0.55cm}@{}}
\toprule
\textbf{Domain} & \textbf{Year} & \textbf{Ref.} & \textbf{What the Hamiltonian represents} & \textbf{Band} & \textbf{Use} \\
\midrule
Audio and acoustics & 2015 & \cite{lopes2015} & Energy of a generic pH audio model; explicit passive integrator & P & S \\
 & 2016 & \cite{falaize2016} & Electrical energy of capacitors, inductors and nonlinear components & P & S \\
 & 2017 & \cite{falaize2017rhodes} & Mechanical and magnetic energy of a Rhodes piano & P & S \\
 & 2020 & \cite{falaize2020speaker} & Electrical, magnetic and suspension energy of a loudspeaker & P & S \\
\midrule
Biology and physiology & 2010 & \cite{jeltsema2010} & Standard pH energy plus memristive ports with no stored energy & P & M \\
 & 2014 & \cite{gawthrop2014} & Chemical free energy of species in biochemical cycles & T & M \\
 & 2015 & \cite{gawthrop2015} & Chemical free energy of hierarchical modules & T & M \\
 & 2016 & \cite{gawthrop2016iet} & Chemical free energy of modular biomolecular systems & T & M \\
 & 2018 & \cite{pan2018cardiac} & Chemical and electrical energy of the cardiac action potential & T & M \\
 & 2021 & \cite{gawthrop2021mod} & Gibbs free energy, constrained by kinetics and data & T & M \\
 & 2022 & \cite{gawthrop2022adv} & Chemical free energy plus hydraulic energy of the carrier fluid & T/P & M \\
\midrule
Chemical reaction networks & 2010 & \cite{vdsMaschke2010crn} & Free energy of an open reaction network & T & M \\
 & 2013 & \cite{vdsRaoJay2013} & Gibbs free energy \eqref{eq:gibbs}, chemical potentials as efforts & T & M \\
 & 2014 & \cite{seslija2014} & Gibbs free energy over spatial compartments & T & M, S \\
 & 2016 & \cite{vdsRaoJay2016} & Gibbs free energy, network dynamics view & T & M \\
 & 2016 & \cite{wang2016fosbe} & Energy and entropy of non-isothermal networks & T & M \\
 & 2018 & \cite{wang2018jmc} & Internal energy and composition, with entropy production & T & M \\
\midrule
Ecology & 2020 & \cite{massaroli2020lv} & Lotka-Volterra first integral \eqref{eq:lv} & A & C \\
\midrule
Economics and markets & 2014 & \cite{macchelli2014econ} & Economic stock; effort times flow is a monetary rate & A & M \\
 & 2017 & \cite{stegink2017tac} & Grid energy plus quadratic storage in prices and multipliers & P/A & C, M \\
\midrule
Fluids, waves and poroelasticity & 2017 & \cite{altmannSchulze2017} & Total energy of a reacting compressible flow & P & M \\
 & 2017 & \cite{cardoso2017} & Kinetic and potential energy of sloshing liquid & P & M \\
 & 2017 & \cite{barradas2017og} & Mechanical and hydraulic energy of a modular wave-energy power take-off & P & M \\
 & 2018 & \cite{almuzakki2018cummins} & Floater energy plus a storage function for wave radiation forces & P & M \\
 & 2021 & \cite{bansal2021} & Kinetic and pressure energy of two phases & P & M \\
 & 2021 & \cite{rashad2021fluid} & Kinetic energy of ideal fluid flow & P & M \\
 & 2021 & \cite{altmann2021poro} & Elastic strain energy plus fluid pressure energy & P & S \\
\midrule
Graphs, crowds and traffic & 2013 & \cite{vdsMaschke2013graphs} & Quadratic disagreement plus edge potentials & A & M, C \\
 & 2013 & \cite{wei2013scl} & Storage (inventory) at network nodes & A & C \\
 & 2023 & \cite{tordeux2023} & Interaction potentials of pedestrians & A & M \\
 & 2024 & \cite{rudiger2024} & Headway potentials plus relative-speed term & A & M \\
 & 2024 & \cite{ackermann2024} & Same, with control ports under noise & A & C \\
\midrule
Infrastructure networks & 2020 & \cite{hauschild2020} & Internal and flow energy of heating networks, with entropy & T & M, S \\
 & 2021 & \cite{strehle2021} & Hydraulic energy of pipe flows & P & M \\
 & 2019 & \cite{phillips2019pumped} & Electrical, mechanical and hydraulic energy of solar pumped-hydro storage & P & M \\
 & 2021 & \cite{krishna2021} & Electrical plus thermal energy of a microgrid & P & C \\
 & 2022 & \cite{machado2022} & Quadratic storage in pipe and tank temperatures & T/A & M \\
 & 2022 & \cite{perryman2022} & Incremental hydraulic energy of a water network & P & C \\
 & 2022 & \cite{esquivel2022pumpturbine} & Electrical, mechanical and hydraulic energy of a reversible pump-turbine & P & M, C \\
\midrule
Irreversible thermodynamics & 2013 & \cite{ramirez2013ces} & Internal energy, with entropy production in the structure & P/T & M \\
 & 2013 & \cite{ramirez2013ejc} & Coupled energies of multi-energy systems & P/T & M, C \\
 & 2013 & \cite{ramirez2013ifac} & Availability as closed-loop storage & T & C \\
 & 2014 & \cite{delvenne2014} & Energy of a pH system in contact with a heat bath & T & M \\
 & 2016 & \cite{ramirez2016aut} & Availability relative to a reference state & T & C \\
 & 2018 & \cite{ramirez2018hyst} & Total energy including heat of a hysteretic actuator & T & M \\
 & 2021 & \cite{lohmayer2021} & Exergy \eqref{eq:exergy} & T & M \\
 & 2021 & \cite{vds2021csm} & Internal energy and entropy on contact geometry & T & M \\
 & 2022 & \cite{ramirez2022ces} & Distributed internal energy with boundary ports & T & M, C \\
 & 2024 & \cite{ramirez2024aut} & Internal energy of interconnected subsystems & T & M \\
\midrule
Machine learning & 2019 & \cite{massaroli2019cdc} & Training objective; network weights are the state & A & L \\
 & 2021 & \cite{desai2021} & Learned Hamiltonian with learned forcing and damping & A & L \\
 & 2022 & \cite{beckers2022cdc} & Gaussian-process distribution over Hamiltonians & P & L \\
 & 2023 & \cite{neary2023} & Learned subsystem Hamiltonians for composition & P & L \\
 & 2026 & \cite{moradi2026} & Learned Hamiltonian under output-error noise & A & L \\
\midrule
Optimization & 2015 & \cite{stegink2015ifac} & Convex objective of the gradient method & A & C \\
 & 2023 & \cite{camlibel2023} & Convex functions of maximal monotone relations & A & M \\
 & 2025 & \cite{gernandt2025} & Primal-dual storage around the optimality system & A & C \\
\midrule
Plasma physics & 2014 & \cite{vu2014ifac} & Magnetic and thermal energy of plasma profiles & P & C \\
 & 2016 & \cite{vu2016mcmds} & Electromagnetic, kinetic and internal energy of the plasma & P & M \\
 & 2016 & \cite{vu2016cep} & Magnetic energy of the safety-factor profile & P & C \\
\midrule
Quantum systems & 2023 & \cite{verrier2023} & Energy functional of the wave function & P & S \\
\midrule
Space instruments & 2024 & \cite{ma2024gw} & Energy of the coupled spacecraft and test-mass dynamics & P & C \\
\bottomrule
\end{tabularx}
\end{table*}

\subsection{A teaching shortlist}
For a graduate course on systems engineering for autonomous systems, a progression that moves from thermodynamic energy to abstract storage and back to physics works well.
\begin{enumerate}[leftmargin=*]
\item \textit{Reaction networks} \cite{vdsRaoJay2016}: free energy replaces mechanical energy without becoming arbitrary.
\item \textit{Irreversible CSTR} \cite{ramirez2013ces, ramirez2016aut}: energy conservation and entropy production are separate statements, and control needs availability rather than energy.
\item \textit{Graphs} \cite{vdsMaschke2013graphs}: one template for consensus and coordination of multi-agent systems.
\item \textit{Predator-prey} \cite{massaroli2020lv}: a non-physical first integral that can still be energy-shaped.
\item \textit{Gradient method and markets} \cite{stegink2015ifac, stegink2017tac}: an algorithm becomes a subsystem that plugs into a plant.
\item \textit{Car following} \cite{rudiger2024}: an abstract Hamiltonian for vehicle autonomy.
\item \textit{Learning} \cite{desai2021, massaroli2019cdc}: pH as a prior, and pH as the training dynamics.
\item \textit{Audio} \cite{falaize2016}: passivity as a numerical stability principle, with every element traceable on a schematic.
\item \textit{Tokamak} \cite{vu2016cep}: IDA-PBC on a reduced distributed model, validated on real hardware.
\end{enumerate}

\section{OPEN QUESTIONS AND CRAZY IDEAS}
\label{sec:Crazy}

The gaps in the literature are as informative as the results.
The items below are search findings, not impossibility claims.
Some are sober extrapolations; others are deliberately speculative.
Each one names the closest verified anchor.

\subsection{One Hamiltonian for compartments}
Three functions in this review look different but are closely related.
Let $\phi(x) = x \ln x - x$ and let $D_\phi(x, y) = \phi(x) - \phi(y) - \phi'(y)(x - y)$ be its Bregman divergence.
The Gibbs free energy \eqref{eq:gibbs} of a reaction network is $\sum_i D_\phi(x_i, x_i^*)$.
The Lotka-Volterra first integral of \eqref{eq:lv}, written in the original coordinates, is $d\, D_\phi(x_1^*, x_1) + b\, D_\phi(x_2^*, x_2)$ up to a constant, with $x_1^* = c/d$ and $x_2^* = a/b$.
The classical Volterra-type Lyapunov functions used for endemic equilibria in mathematical epidemiology have the same form.
So the same relative-entropy divergence, with its arguments in one order or the other, serves as the Hamiltonian for chemistry, ecology and epidemics.

\begin{conjecture}[Epidemics as open pH systems]
An SIR or SEIR model with vital dynamics is a mass-action network with irreversible reactions ($S + I \to 2I$, $I \to R$). It should admit a pH or irreversible-pH representation with a relative-entropy Hamiltonian, in which vaccination, testing and travel enter as ports. Endemic-level targeting would then become an energy-shaping problem.
\end{conjecture}
The obstacle is real.
Infection is irreversible, so the network is not detailed-balanced and \eqref{eq:crn} does not apply directly.
The irreversible pH machinery \cite{ramirez2013ces, ramirez2022entropy} and the extinction-free control of \cite{massaroli2020lv} are the obvious starting points.
We found Hamiltonian reformulations of epidemic models only as preprints, without ports and without control.

\subsection{Supply chains with prices as efforts}
Inventories are conserved stocks, shipments are edge flows and prices are natural dual variables.
The graph template \cite{vdsMaschke2013graphs}, the load-balancing controllers of \cite{wei2013scl}, the market coupling of \cite{stegink2017tac} and the monotone theory of \cite{camlibel2023} together suggest that a production network and its pricing mechanism could be one passive interconnection, with a single stability proof.
Macchelli's macroeconomic model \cite{macchelli2014econ} shows that monetary port pairs are workable.
We found no peer-reviewed pH supply-chain model.

\subsection{Opinion dynamics with stubborn agents}
Consensus is covered \cite{vdsMaschke2013graphs}; opinion formation is not.
In Friedkin-Johnsen-type models each agent is anchored to a private prejudice.
In pH language this anchoring looks like a resistor connected to a constant source, and a news outlet looks like an input port.
Polarisation would then be an equilibrium of a shaped storage function, and the question of whether a moderating input can reshape it becomes a control problem.
This is speculative, but the structure is suggestive.

\subsection{Optimizers as mechanical systems}
The continuous-time heavy-ball method, $\ddot{z} + \gamma \dot{z} + \nabla f(z) = 0$, is a damped mechanical system.
With $p = \dot{z}$ and $H = \tfrac{1}{2}\norm{p}^2 + f(z)$ it is exactly \eqref{eq:sigma_total} with canonical $J$ and $R = \mathrm{diag}(0, \gamma I)$, so the objective is literally potential energy.
Primal-dual methods are already pH \cite{gernandt2025}, and so is one view of network training \cite{massaroli2019cdc}.
An open question is which widely used optimisers with adaptive step sizes admit useful pH realisations, and whether the structure brings robustness guarantees rather than only new notation.

\subsection{Learned pH digital twins for structural health monitoring}
A learned pH model separates storage, interconnection and dissipation by construction \cite{desai2021, neary2023, moradi2026}.
In a structure, damage changes stiffness, which lives in $H$, or damping, which lives in $R$.
A pH network identified on healthy data and then re-identified on monitoring data could therefore localise and classify damage by tracking which part of the model drifts.
Gaussian-process pH models would add calibrated uncertainty to that decision \cite{beckers2022cdc}, and compositional learning would let each structural member keep its own learned Hamiltonian \cite{neary2023}.
A model-based first step exists.
For discrete Timoshenko beams, the violation of the power balance \eqref{eq:power_balance} serves as a consistency residual that splits into element-level damage signatures and locates stiffness loss without prior knowledge of where it occurred \cite{esquivel2026timoshenko}.
A learned version, identified from monitoring data rather than from a finite-element model, has not yet been tried.

\subsection{Further gaps}
Open quantum systems need a pH theory for Lindblad dissipation in which complete positivity, not just passivity, is preserved; the closed Schrödinger case exists \cite{verrier2023}.
Whole-body cardiovascular models have genuine pressure-flow power ports and compliant storage, yet only cell-level bond-graph work was found \cite{pan2018cardiac}.
Macroscopic traffic models of Lighthill-Whitham-Richards type with tolling or routing optimisers look like candidates for monotone pH interconnection \cite{camlibel2023, rudiger2024}.
Finally, an autonomous robot whose storage function is its remaining exergy \cite{lohmayer2021}, combining battery, thermal and mechanical budgets in one Hamiltonian, would turn energy-aware autonomy into a passivity problem.

\section{CONCLUDING REMARKS}
\label{sec:Concluding_Remarks}

The port-Hamiltonian framework has moved along a continuum, from literal energy, through thermodynamic potentials, to storage functions with no physical meaning.
Chemistry and irreversible thermodynamics are where the theory has evolved most, because they forced a rethink of what the Hamiltonian should be.
Plasma, audio, wave energy and drag-free satellites show that pH structure can organise very complicated real-energy models, and plasma control has run on real hardware.
Economics, ecology, optimisation and learning are where it has been stretched furthest.

What survives at the abstract end is not energy but structure: modular interconnection through conjugate ports, a dissipation inequality, a clean separation of conservative and dissipative dynamics, and stability proofs that survive interconnection.
The monotonicity results of \cite{camlibel2023} give this abstract view a solid mathematical foundation.

The most promising open territory is where a compartmental or network structure meets a convex potential: epidemics, supply chains, opinion dynamics and traffic.
In those domains the useful question is not whether the equations can be rewritten in Hamiltonian form.
It is whether a pH storage and port structure exposes a conservation law, a dissipation inequality or a control principle that was not already obvious.

\bibliographystyle{IEEEtran}
\bibliography{refs}


\end{document}